\documentclass[journal]{IEEEtran}

\usepackage{cite}
\usepackage{color}

\ifCLASSINFOpdf
   \usepackage[pdftex]{graphicx}
\else
\fi

\usepackage{amsmath}

\begin{document}

\title{Distribution Network Marginal Costs: Enhanced AC OPF Including Transformer Degradation}

\author{Panagiotis~Andrianesis,~\IEEEmembership{Member,~IEEE},
        and~Michael~Caramanis,~\IEEEmembership{Senior~Member,~IEEE}
\thanks{P. Andrianesis and M. Caramanis are with the Systems Eng. Div., Boston University, Boston, MA: panosa@bu.edu, mcaraman@bu.edu. Research partially supported by Sloan Foundation G-2017-9723 and NSF AitF 1733827.
}
}

\maketitle
\begin{abstract}
This paper considers the day-ahead operational planning problem of radial distribution networks hosting Distributed Energy Resources, such as Solar Photovoltaic (PV) and Electric Vehicles (EVs). 
We present an enhanced AC Optimal Power Flow (OPF) model that estimates dynamic Distribution nodal Location Marginal Costs (DLMCs) encompassing transformer degradation as a short-run network variable cost.
We decompose real/reactive power DLMCs into additive marginal cost components representing respectively the costs of real/reactive power transactions at the T\&D interface, real/reactive power marginal losses, voltage and ampacity congestion, and transformer degradation, capturing intertemporal coupling. Decomposition is useful in identifying sources of marginal costs and facilitating the employment of distributed AC OPF algorithms. DLMCs convey sufficient information to represent the benefit of shifting real/reactive power across time and achieve optimal Distribution Network and DER operation.
We present and analyze actual distribution feeder based numerical results involving a wealth of future EV/PV adoption scenarios and illustrate the superiority of our approach relative to reasonable conventional scheduling  alternatives.
Overwhelming evidence from extensive numerical results supports the significant benefits of internalizing short-run marginal asset ---primarily transformer--- degradation.
\end{abstract}

\begin{IEEEkeywords}
Distribution Locational Marginal Costs, Distributed Energy Resources, Marginal Transformer Degradation.
\end{IEEEkeywords}

\IEEEpeerreviewmaketitle

\section{Introduction}

\IEEEPARstart{R}{apidly} growing adoption of Distributed Energy Resources (DERs), including clean, albeit volatile, renewable generation, 
storage and flexible loads with storage-like properties and Volt/VAR control capabilities, e.g., Electric Vehicles (EVs) and Solar Photovoltaic (PV) inverters, present a major challenge together with still unexploited opportunities.

Numerous studies on DER integration \cite{OlivierEtAl_2016, GongEtAl_2012, HilsheyEtAl_2013, QianEtAl_2015, AssolamiMorsi2015, ElNozahySalama_2014, PezeshkiEtAl_2014, GrayMorsi_2017, WangEtAl_2018} address the impact of various DERs (primarily EVs and PVs) on the Distribution Network and its assets, e.g., overvoltages due to PV, and explore DER capabilities, e.g., the provision of reactive power. 
Interestingly, several works investigate the impact of increasing EV penetration on distribution transformers \cite{GongEtAl_2012, HilsheyEtAl_2013, QianEtAl_2015, AssolamiMorsi2015, ElNozahySalama_2014}, noting that transformer aging is dependent upon the thermal effects of persistent loading. 
Indeed, the transformer life is strongly related to its winding hottest spot temperature (HST) driving insulation aging; assuming operation at a reference HST, it may exceed 20 years.
IEEE C.57.91-2011 \cite{IEEE2011} and IEC 60076-7:2018 \cite{IEC2018} standards provide an exponential representation of the aging acceleration factor, which is greater (less) than 1 for HSTs exceeding (being lower than) $110^o$C. 
Early EV pilot studies \cite{EVproject2013} identified that clustering EV chargers under the same transformer may cause damage from persistent overloading.
A recent study for the Sacramento Municipal Utility District \cite{SMUD} estimated that
about 26\% of the substations and service transformers would experience overvoltages  
due to PVs, whereas up to 17\% of the approximately 12,000 service transformers would experience overloads  
due to EVs and would need to be replaced at an average estimated cost of \$7,400 per transformer.

There is also an emerging literature that extends wholesale Locational Marginal Price (LMP) markets to distribution networks through reliance on Distribution LMPs (DLMPs).
This paper argues that most published work is lacking in considering the daily cycle costs of distribution networks and complex DER preferences that are inadequately represented by uniform price quantity bids/offers.
Indeed, past work has paid insufficient attention to Grid-DER co-optimization and the role of economically efficient Distribution Locational Marginal Costs (DLMCs) that go beyond wholesale market LMPs by internalizing \emph{intertemporally} coupled DER preferences, e.g., energy rather than capacity demanding EVs, and salient distribution network costs, e.g., transformer Loss of Life (LoL).
The considerable time elapsed between the proposal of wholesale power markets in the 80s \cite{CaramanisEtAl_1982, SchweppeEtAl_1988} and their eventual extensive and successful implementation in the late 90s, behooves us to redouble our effort in the appropriate modeling and estimation of DLMCs.
The expected shock that distributed generation and mobile EV battery charging will most likely deliver to distribution assets (transformers in particular) may dwarf the shock delivered in the late 70s by the massive adoption of air-conditioning.
What costs should DLMCs capture to promote efficient scheduling of EV charging and dual use of PV smart inverters?

Existing literature investigating the EV impact on distribution transformers \cite{GongEtAl_2012, HilsheyEtAl_2013, QianEtAl_2015, AssolamiMorsi2015, ElNozahySalama_2014} includes mostly simulation studies assessing the impact of various charging schedules on transformer LoL.
Indicative results show that the daily LoL of a residential service transformer may almost double if EVs charge upon arrival \cite{QianEtAl_2015}, and that ``open-loop'' EV charging, e.g., delay until after midnight that may be promoted by Time of Use rates, can actually increase, rather than decrease, transformer aging \cite{HilsheyEtAl_2013}. 
Rooftop PVs, on the other hand, can utilize their inverters to mitigate overvoltage issues and most importantly to schedule their real and reactive power output with a beneficial impact on transformer LoL.
Indeed, \cite{PezeshkiEtAl_2014} shows that high PV penetration can significantly extend the life of distribution transformers in a suburban area.
Furthermore, while previous research determined insignificant synergy between the substation transformer LoL and charging EVs in the presence of rooftop PV (due to non-coincidence between peak hours of PV generation and EV charging), \cite{GrayMorsi_2017} shows that the transformer thermal time constant allows PV generation to reduce the  HST when EVs are charging.

Although the above works consider the DER impact on transformer aging, they do not internalize degradation to the operational planning problem.
A first outline is provided in \cite{CaramanisEtAl2016, TaborsEtAl2016, NtakouEtAlIEEETD2014}, among the first DLMP works in which a transformer model is employed in AC Optimal Power Flow (OPF) with nonlinear LoL objective function components and omission of the thermal time constant intertemporal coupling.
These works employ the \emph{branch} \emph{flow} (a.k.a. DistFlow) equations, introduced for radial networks by \cite{BaranWu1989}, and their relaxation to convex Second Order Cone constraints proposed by \cite{FarivarLow2013}.
Parallel works \cite{LiEtAl_2014, HuangEtAl_2015} determined DLMPs 
using DC OPF, considering price-taking EV aggregators.
Quadratic programming is used in \cite{HuangEtAl_2015}, as opposed to linear programming in \cite{LiEtAl_2014}, to derive DLMPs that are announced to EV and heat pump aggregators optimizing their energy plans. 
DLMP approaches using linearized models are discussed in \cite{YuanEtAl_2018}, whereas \cite{HanifEtAl_2019} also employs a trust-region method. 
Following \cite{CaramanisEtAl2016}, \cite{BaiEtAl-DLMP}  employs the relaxed DistFlow model and proposes a day-ahead framework, in which DERs bid into a distribution-level electricity market.
It acknowledges that reactive power and voltage are critical features in distribution network operation, and, hence, DC OPF approaches are inadequate in capturing significant costs.
DLMPs are derived using sensitivity factors of the linear version of the DistFlow model;
they include ampacity congestion but do not consider the impact on transformers.
In a parallel to \cite{BaiEtAl-DLMP}, \cite{Papavasiliou_2017} provides a comprehensive analysis of various approaches in decomposing and interpreting DLMPs, and discusses accuracy and computational aspects.

Our main objective is to derive the time and location specific marginal costs of real and reactive power in the distribution network that are consistent with the optimal DER schedule in the day-ahead operational planning problem. 
We aim at discovering DLMCs and understanding the components, building blocks, and sources that constitute the marginal costs. 
Even more importantly, this decomposition enables scaling to real distribution feeders by leveraging distributed algorithms requiring a simple day-ahead AC OPF solution conditional upon given DER schedules to derive tentative DLMCs that drive individual daily operational DER scheduling. 
The main contributions are three-fold.

First, we present an enhanced AC OPF model that includes a detailed cost representation of transformer degradation in the presence of DERs (PVs, EVs).
More specifically, we enhance the relaxed branch flow model  with accurate modeling of the transformer LoL and intertemporal thermal response.
Second, we provide intuitive formulas of the real and reactive power DLMC component representing transformer degradation cost.
Employing sensitivity analysis, we decompose DLMCs into additive cost components, based on the cost of real and reactive power at the substation, real and reactive power marginal losses, voltage congestion, line (ampacity) congestion, and transformer degradation.
We illustrate how the last component captures the intertemporal impact of the thermal response dynamics on the marginal cost and quantify this impact by showing how the transformer load at a specific hour affects future hour degradation marginal costs.
We also discuss practical considerations for capturing the transformer degradation impact that extends beyond the day-ahead horizon. 
Third, we employ an actual distribution feeder and a wealth of future DER adoption scenarios, 
and we compare to approaches that are often considered reasonable (Business as Usual, 
time varying LMPs that do not vary across distribution nodes, and traditional Line Loss minimization).
Extensive numerical results illustrate the significant benefits of internalizing short-run marginal asset, primarily transformer, degradation, and that voltage and ampacity congestion DLMC components are insufficient in providing the price signal that supports the system optimal schedules.
They also demonstrate that apart from optimal short-run scheduling, our approach can harvest otherwise idle DER reactive power compensation capabilities, and increase distribution network DER (EV/PV) hosting capacity mitigating investments in distribution infrastructure.

The remainder of this paper is organized as follows.
Section \ref{Model} introduces the enhanced AC OPF model.
Section \ref{DLMCs} presents the decomposition of DLMCs. 
Section \ref{CaseStudy} discusses the case study based numerical findings.
Section \ref{Conclusions} concludes and provides further research directions.

\section{Enhanced AC OPF Model} \label{Model}

In this section, we introduce the network and DER models (Subsection \ref{NetworkModel}), we present the detailed transformer degradation formulation 
(Subsection \ref{TransformerModel}), and we summarize the enhanced AC OPF problem (Subsection \ref{ProblemSummary}) for clarity.

\subsection{Network and DER Models} \label{NetworkModel}

We consider a radial network with $N+1$ nodes and $N$ lines.
Let $\mathcal{N}= \{ 0,1,...,N \}$ be the set of nodes, with node $0$ representing the root node, and $\mathcal{N^+} \equiv \mathcal{N} \setminus \{0\}$.
Let $\mathcal{L}$ be the set of lines, with each line denoted by the pair of nodes $(i,j)$ it connects --- henceforth $ij$ for short, where node $i$ refers to the (unique due to the radial structure) preceding node of $j \in \mathcal{N^+}$.
Transformers are represented as a subset of lines, denoted by $y \in \mathcal{Y} \subset \mathcal{L}$. 
For node $i \in \mathcal{N}$, $v_i$ denotes the magnitude squared voltage.
For node $j \in \mathcal{N^+}$, $p_j$ and $q_j$ denote the net demand of real and reactive power, respectively.
A positive (negative) value of $p_j$ refers to withdrawal (injection); 
similarly for $q_j$.
Net injections at the root node are denoted by $P_0$ and $Q_0$, for real and reactive power, respectively.
These are positive (negative) when power is flowing from (to) the transmission system.
For each line $ij$, with resistance $r_{ij}$ and reactance $x_{ij}$, $l_{ij}$ denotes the magnitude squared current, $P_{ij}$ and $Q_{ij}$ the sending-end real and reactive power flow, respectively.

The \emph{branch} \emph{flow} (AC power flow) 
equations, are listed below, where we introduce the time index $t$;
unless otherwise mentioned, $j \in \mathcal{N^+}$, and $t \in \mathcal{T^+}$, with $\mathcal{T} = \{0,1,...,T \}$, $\mathcal{T^+} \equiv \mathcal{T} \setminus \{ 0 \}$, and $T$ the length of the optimization horizon.
\begin{equation} \label{EqSubstation}
P_{0,t} = P_{01,t}  \rightarrow (\lambda_{0,t}^P), \quad Q_{0,t} = Q_{01,t}  \rightarrow(\lambda_{0,t}^Q), \quad \forall t,
\end{equation}
\begin{equation} \label{EqRealBalance}
P_{ij,t} - r_{ij} l_{ij,t} = \sum_{k: j \to k} P_{jk,t} + p_{j,t} \rightarrow (\lambda_{j,t}^P), \quad \forall j,t,
\end{equation}
\begin{equation} \label{EqReactiveBalance}
Q_{ij,t} - x_{ij} l_{ij,t} = \sum_{k: j \to k} Q_{jk,t} + q_{j,t} \rightarrow (\lambda_{j,t}^Q), \quad \forall j,t,
\end{equation}
\begin{equation} \label{EqVoltageDef}
v_{j,t} = v_{i,t} - 2 r_{ij} P_{ij,t} - 2 x_{ij} Q_{ij,t} + \left( r_{ij}^2 + x_{ij}^2 \right) l_{ij,t}, \,\,\, \forall j,t,
\end{equation}
\begin{equation} \label{EqCurrentDef}
 v_{i,t} l_{ij,t} = P_{ij,t}^2 + Q_{ij,t}^2 \quad \forall j,t.
\end{equation}
Briefly, (\ref{EqSubstation})--(\ref{EqReactiveBalance}) define the real and reactive power balance, (\ref{EqVoltageDef}) the voltage drop, and (\ref{EqCurrentDef}) the apparent power but can be also viewed as the definition of the current.
The real (reactive) power net demand $p_{j,t}$ ($q_{j,t}$) includes the aggregate effect of: 
(\emph{i}) conventional demand consumption $p_{d,t}$ ($q_{d,t}$) of load $d \in \mathcal{D}_j$, where $\mathcal{D}_j \subset \mathcal{D} $ is the subset of loads (set $ \mathcal{D}$) connected at node $j$;
(\emph{ii}) consumption $p_{e,t}$ ($q_{e,t}$) of EV $e \in \mathcal{E}_{j,t}$, where $\mathcal{E}_{j,t} \subset \mathcal{E}$ is the subset of EVs (set $\mathcal{E}$) that are connected at node $j$, during time period $t$, and 
(\emph{iii}) generation $p_{s,t}$ ($q_{s,t}$) of PV (rooftop solar) $s \in \mathcal{S}_j$, where $\mathcal{S}_j \subset \mathcal{S}$ is the subset of PVs (set $\mathcal{S}$) connected at node $j$.
For clarity, the definitions of aggregate dependent variables are listed below:
\begin{equation} \label{Pinject}
p_{j,t} = \sum_{d \in \mathcal{D}_j} p_{d,t} + \sum_{e \in \mathcal{E}_{j,t}}  p_{e,t} - \sum_{s \in \mathcal{S}_j} p_{s,t}, \quad\forall j, t,
\end{equation}
\begin{equation} \label{Qinject}
q_{j,t} = \sum_{d \in \mathcal{D}_j} q_{d,t} + \sum_{e \in \mathcal{E}_{j,t}}  q_{e,t} - \sum_{s \in \mathcal{S}_j} q_{s,t} , \quad \forall j, t.
\end{equation}

We also supplement with voltage and current limits: 
\begin{equation} \label{EqVoltageLimits}
\underline{v}_i \leq v_{i,t} \leq \bar{v}_i \rightarrow ( \underline{\mu}_{i,t}, \bar \mu_{i,t}), \quad \forall i,t,
\end{equation}
\begin{equation} \label{EqCurrentLimits}
l_{ij,t} \leq \bar{l}_{ij}  \rightarrow ( \bar{\nu}_{j,t} ), \quad  \forall j,t,
\end{equation}
where $\underline{v}_i$, $\bar{v}_i$, and $\bar{l}_{ij}$ are the lower voltage, upper voltage, and line ampacity limits (squared), respectively.
Dual variables of constraints (\ref{EqSubstation})--(\ref{EqReactiveBalance}), (\ref{EqVoltageLimits}) and (\ref{EqCurrentLimits}) are shown in parentheses.

Lastly, we provide PV and EV constraints for a general multi-period operational planning problem setting, accommodating smart inverter capabilities and EV mobility.

Consider a PV with nameplate capacity $C_s$, irradiation level $\rho_t \in [0,1]$, and $\mathcal{T}_I \subset \mathcal{T^+}$ the subset of time periods for which $\rho_t > 0$.
PV constraints ($\forall s \in \mathcal{S}$) are as follows:
 \begin{equation} \label{EqPVcon1}
 p_{s,t} \leq \rho_t C_{s}, \quad p_{s,t}^2 + q_{s,t}^2 \leq C_s^2, \quad \forall s, t \in \mathcal{T}_I,
 \end{equation}
\begin{equation} \label{EqPVcon2}
p_{s,t} = q_{s,t} = 0, \qquad \forall s, t \not \in \mathcal{T}_I,
\end{equation}
with $p_{s,t} \ge 0$. 
Constraints (\ref{EqPVcon1}) impose limits on real and apparent power (implicitly assuming an appropriately sized inverter), whereas (\ref{EqPVcon2}) imposes zero generation when $\rho_t = 0$.

Consider an EV connected for $Z$ intervals, at nodes $j_1,...,j_Z$.
Let $\mathcal{T}^{beg} = \{ \tau_z^{beg} \}$ ($\mathcal{T}^{end} = \{ \tau_z^{end} \}$) be the set of time periods of interval $z$, for $z=1,...,Z$, denoting an adjusted beginning (end), considering the part of the interval within the time horizon. 
Let $\mathcal{T}_z = \{ \tau_z^{beg}+1,...,\tau_z^{end} \}$ be the set of time periods of interval $z$, during which the EV is connected at node $j_z$.
The State of Charge (SoC), $u_{e,t}$, for EV $e$, $t \in \mathcal{T}_e^{beg} \cup \mathcal{T}_e^{end}$ (index $e$ was omitted earlier for simplicity), is reduced by $\Delta u_{z,e}$ after departure $z$ and until arrival $z+1$, for $z = 1,..,Z_e-1$. 
EV constraints ($\forall e \in \mathcal{E}$) are as follows:
\begin{equation} \label{EVCon2}
u_{e,\tau_z^{end}} = u_{e,\tau_z^{beg}} + \sum_{t \in \mathcal{T}_{e,z}} p_{e,t}, \quad \forall e, z = 1,...,Z_e,
\end{equation}
\begin{equation} \label{EVCon3}
u_{e,\tau_{z+1}^{beg}} = u_{e,\tau_z^{end}} - \Delta u_{e,z}, \quad \forall e, z = 1,...,Z_e-1,
\end{equation}
\begin{equation} \label{EVCon4} 
u_{e,\tau_1^{beg}} = u_e^{init}, \quad u_{e,t}^{min} \leq u_{e,t} \leq C_e^B, \quad \forall e, t \in \mathcal{T}_{e}^{end},
\end{equation}
\begin{equation} \label{EVCon5}
p_{e,t}^2 + q_{e,t}^2 \leq C_e^2, \quad 0 \leq p_{e,t} \leq C_{r}, \quad \forall e, t \in { \cup_{z=1}^{Z_e} \mathcal{T}_{e,z} },
\end{equation}
\begin{equation} \label{EVCon6}
 p_{e,t} = q_{e,t} = 0, \,\forall e, t \in {\mathcal{T^+}} \setminus { \cup_{z=1}^{Z_e} \mathcal{T}_{e,z} },
\end{equation}
with $p_{e,t}, u_{e,t} \geq 0$.
Eqs. \eqref{EVCon2} and \eqref{EVCon3} define the SoC at the end/beginning of an interval, after charging/traveling, respectively. 
Constraints \eqref{EVCon4} initialize the SoC ($u_e^{init}$) at $\tau_1^{beg}$ and impose a minimum SoC, $u_{e,t}^{min}$, at the end of an interval as well as the limit of the EV battery capacity, $C_e^B$. 
Constraints \eqref{EVCon5} impose the limits of the charger, $C_e$, (related to the size of the inverter) and the charging rate, $C_{r}$, (related to the capacity of the EV battery charger).
Lastly, \eqref{EVCon6} imposes zero consumption when the EV is not plugged in.

\subsection{Transformer Degradation Formulation} \label{TransformerModel}

IEEE and IEC Standards \cite{IEEE2011, IEC2018} provide detailed formulas on the transformer thermal response, widely used in simulation studies, but not employed in an AC OPF model.
For a given HST of the winding, $\theta^H$, they define the aging acceleration factor as $F_{AA} = \exp \big( \frac{15000}{383}-\frac{15000}{\theta^H + 273} \big)$.
We consider a piecewise linear approximation, $\tilde{F}_{AA}$, 
given by:
\begin{equation} \label{Eq:PieceWise} 
\tilde{F}_{AA} = a_{\kappa} \theta^H - b_{\kappa}, \,\, \theta_{\kappa-1}^H \leq \theta^H < \theta_{\kappa}^H, \,\, \kappa = 1,...,M,
\end{equation}
where $M$ is the number of segments.
Notably, appropriate selection of the breakpoints ensures a highly accurate representation of the exponential aging model.
Introducing indices $y$ and $t$, and substituting $\tilde{F}_{AA}$ with $f_{y,t}$ to simplify the notation, we replace \eqref{Eq:PieceWise} with the following set of inequalities:
\begin{equation} \label{PiecewiseIneq}
f_{y,t} \geq a_{\kappa} \theta^H_{y,t} - b_{\kappa}, \qquad \forall y,t,\kappa,
\end{equation}
with $f_{y,t} \geq 0$. 
Since $f_{y,t}$ will be included with a cost in the objective function, at least one of these inequalities should be binding, 
and hence, the relaxation of \eqref{Eq:PieceWise} to the linear inequality constraints \eqref{PiecewiseIneq}, which would otherwise require the introduction of binary variables, is exact.

The transformer HST, $\theta^H_{y,t}$, is given by: 
\begin{equation} \label{EqThetaH}
\theta^H_{y,t} = \theta^A_{y,t} + \Delta \theta^{TO}_{y,t} + \Delta \theta^H_{y,t} = \theta^{TO}_{y,t} + \Delta \theta^H_{y,t},
\end{equation}
where
$\theta^A_{y,t}$ is the ambient temperature at the transformer location, $\theta^{TO}_{y,t}$ is the top-oil (TO) temperature,
$\Delta \theta^{TO}_{y,t}$ is the TO temperature rise over $\theta^A_{y,t}$,
and $\Delta \theta^H_{y,t}$ is the winding HST rise over $\theta^{TO}_{y,t}$.
The detailed derivations are based on the heat transfer differential equations, and they are described in \cite{IEC2018} as a function of time, for varying load and ambient temperature, using exponential and difference equations.

For the TO temperature, $\theta^{TO}_t$ (we dropped index $y$), the differential equation is given by:
\begin{equation} \label{EqDif}
\Delta \bar \theta^{TO} \Big(\frac{1 + K_t^2 R}{1 + R}\Big)^n = k_{1} \tau^{TO}\frac{d\theta^{TO}_t}{dt} + \theta^{TO}_t - \theta^A_t,
\end{equation}
where $\Delta \bar \theta^{TO}$ is the rise of TO temperature over ambient temperature at rated load, $K_t$ is the ratio of the (current) load to the rated load, $R$ is the ratio of load losses at rated load to no-load losses, $\tau^{TO}$ is the oil time constant with recommended value 3 hours, $k_{1}$ and $n$ are constants with recommended values 1 and 0.8, respectively.
Since the granularity of our problem is less than half the recommended value of $\tau^{TO}$ ($\Delta t = 1$ hour), employing the difference equations, we get:
\begin{equation} \label{Eq:thetaTO}
\theta_t^{TO} = \delta \theta_{t-1}^{TO} 
  + (1-\delta) \Big[ \Delta \bar \theta^{TO} \Big(  \frac{1+K_t^2 R}{1+R}\Big)^n  + \theta_t^{A} \Big],
\end{equation}
where $\delta = \frac{k_{1} \tau^{TO}}{k_{1} \tau^{TO} + \Delta t}$.
For the winding HST rise over $\theta^{TO}_t$, $\Delta \theta^H_t$, it can be shown that \cite{IEEE2011} and \cite{IEC2018} yield the same results for distribution (small) transformers.
Because the winding time constant $\tau^w$ has an indicative value of about 4 min (much less than the hourly granularity), the transient behavior $1-\exp \left(-\frac{\Delta t}{\tau^w} \right)  \approx 1$, vanishes.
Hence, using \cite{IEEE2011}, we get:
\begin{equation} \label{Eq:DthetaH}
\Delta \theta_t^{H} = \Delta \bar \theta^{H} \big( K_t^2\big) ^{m},
\end{equation}
where $\Delta \bar \theta^{H}$ is the rise of HST over TO temperature at rated load, and $m$ is a constant with recommended value 0.8.

The load ratio $K_{t}$ can be defined \emph{w.r.t.} the transformer nominal current (at rated load), denoted by $I^N$.
Using variables $l_{t}$ (omitting index $y$), we have $K_{t}^2 = {l_{t}}/{l^N}$, where $l^N = ({I^N})^2$.
Hence, approximating the terms $\left(\frac{1+K_t^2 R}{1+R}\right)^n$ and $\left( K_t^2\right) ^{m}$ in (\ref{Eq:thetaTO}) and (\ref{Eq:DthetaH}), respectively, using the 1st order Taylor expansion of $K_{t}^2$ around 1, equivalently of $l_t$ around $l^N$ 
 --- one can select a different point at each time period depending on the anticipated loading conditions,
and replacing $K_{t}^2$ by ${l_{t}}/{l^N}$, we get
$\left(\frac{1+{K_t}^2 R}{1+R} \right)^n \approx  
 \frac{nR}{(1+R)} \frac{l_t}{l^N} + \frac{1+(1-n)R}{1+R}$,
$\left( K_t^2\right) ^{m}  \approx m  \frac{l_t}{l^N} + 1 - m$,
and \eqref{Eq:thetaTO} and \eqref{Eq:DthetaH} yield:
\begin{equation} \label{Eq:thetaTOapprox}
\begin{split}
\theta_t^{TO}  = \delta \theta_{t-1}^{TO}
& + (1-\delta)
    \frac{n R \Delta \bar \theta^{TO}}{1 + R}  \frac{l_t}{l^N} \\ 
& + (1-\delta) 
    \Big[
    \frac{ 1 + (1-n)R }{1+R} \Delta \bar \theta^{TO}
    + \theta_t^{A} \Big],
\end{split} 
\end{equation}
\begin{equation} \label{Eq:DthetaHapprox}
\Delta \theta_t^{H} =  m \Delta \bar \theta^{H} \frac{l_t}{l^N}  + (1-m) \Delta \bar \theta^{H}.
\end{equation}

We highlight two important features that greatly facilitate the accurate representation of the transformer thermal response.
The first feature is related to the hourly granularity of the operational planning problem. 
Both $\theta^{TO}_{y,t}$ and $\Delta \theta^H_{y,t}$ are characterized by thermal time constants, whose typical values (about 3 hours for the oil and about 4 minutes for the winding), allow us to employ difference equations for the oil whereas assume a steady state for the winding.
Notably, this would still hold if we allowed up to a 15-min granularity.
The second feature is related to the magnitude squared current variable, $l_{y,t}$, of the branch flow model, which allows us to define the square of the ratio of the transformer load to the rated load using $l_{y,t}$ and the nominal (at rated load) current squared, denoted by $l_y^N$.
Hence, the equations that define the HST fit nicely with the branch flow model. 

Combining \eqref{EqThetaH} and \eqref{Eq:DthetaHapprox}, we define $\theta_{y,t}^{H}$ as a linear equation with variables $\theta_{y,t}^{TO}$ and $l_{y,t}$. 
Substituting $\theta_{y,t}^{TO}$ with $h_{y,t}$, to simplify the notation, and replacing $\theta_{y,t}^H$ in \eqref{PiecewiseIneq}, we get:
\begin{equation} \label{Xf1}
f_{y,t} \geq \alpha_{\kappa} h_{y,t} + \beta_{y, \kappa} l_{y,t} + \gamma_{y,\kappa}
\rightarrow (\xi_{y,t,\kappa}), \, \forall y, t, \kappa,
\end{equation}
with $f_{y,t} \ge 0$, $\xi_{y,t,\kappa}$ the dual variable, and coefficients $\alpha_\kappa$, $\beta_{y,\kappa}$,  $\gamma_{y,\kappa}$ obtained directly from \eqref{PiecewiseIneq} and \eqref{Eq:DthetaHapprox}, with recommended values $\alpha_\kappa = a_{\kappa}$,
$\beta_{y,\kappa} =  a_{\kappa} \frac{4 \Delta \bar{\theta}_{y}^H}{ 5 l_y^N}$, 
$\gamma_{y,\kappa} = a_{\kappa} \frac{\Delta \bar\theta_{y}^H}{5} - b_{\kappa}$. 

The TO temperature at time period $t$, $\theta_{y,t}^{TO}$, is in turn defined by a linear recursive equation that includes the TO temperature at time period $t - 1$, $\theta_{y,t-1}^{TO}$, and $l_{y,t}$. Indeed, there is also an impact of the ambient temperature, $\theta_{y,t}^{A}$, but this is in fact an input parameter.
Substituting $\theta_{y,t}^{TO}$ with $h_{y,t}$ in \eqref{Eq:thetaTOapprox}, we get:
\begin{equation} \label{Xf2}
h_{y,t} = \delta h_{y,t-1} + \epsilon_y l_{y,t} + \zeta_{y,t},\quad \forall y, t,
\end{equation}
where the coefficients $\delta$, $\epsilon_{y}$, and $\zeta_{y,t}$ are obtained directly from \eqref{Eq:thetaTOapprox}, with recommended values $\delta = \frac{3}{4}$,
$\epsilon_{y} =  \frac{ R_y \Delta \bar \theta_{y}^{TO} }{5(1 + R_y) l_y^N }$,
$\zeta_{y,t} = \frac{(5 + R_y)\Delta \bar \theta_{y}^{TO}}{20(1 + R_y)} + \frac{\theta_{y,t}^{A}}{4}$. 
Notably, the value of $\delta$ that is always less than $1$ (for $\Delta t$ equal to 1 h, 30 min or 15 min, $\delta$ will be 0.75, 0.857 or 0.923, respectively) plays an important role in the transformer degradation cost intertemporal impact.  

One issue naturally arising in the context of the enhanced operational planning problem that internalizes intertemporal transformer cost relations is the impact of decisions at period $t$ extending beyond the optimization horizon.
First, rolling horizon approaches are typical 
potential remedies,
which would also provide an additional means for dealing with day-ahead to real-time uncertainty associated with day-ahead forecasted and actual real-time DER preferences.
Second, an extended horizon for the transformer degradation cost might offer a reasonable alternative solution; it should include the next half day to capture 4 oil time constants, while accounting for anticipated ambient temperature trend.
A third approach, quite suitable for simulations, is to set the following constraint:
\begin{equation} \label{Cycle}
h_{y,T} = h_{y,0} \to (\rho_y), \qquad \forall y, 
\end{equation}
which essentially models the daily 24-hour ahead problem as a cycle repeating over identical days.
This constraint's dual variable, $\rho_y$, captures the future impact on  transformer LoL of loading towards the end of the day. 
Alternatively, we could require some explicit condition for $h_{y,T}$, e.g., equal or less than or equal to some target for the next day initial condition, and add a cost in the objective function (pre-calculated by offline studies) for deviating from this target.

\subsection{Optimization Problem (Enhanced AC OPF) Summary} \label{ProblemSummary}

The objective function of the enhanced AC OPF problem --- referred to as \textbf{Full-opt}--- aims at minimizing the aggregate real and reactive power cost, and the transformer degradation cost. 
The real power cost, $c_t^P$, is typically the LMP at the T\&D interface, $c_t^Q$ can be viewed as the opportunity cost for the provision of reactive power, and $c_y$ is the hourly cost for transformer $y$.
Full-opt is defined as follows:
\begin{equation} \label{Obj2}
\text{minimize}
\overbrace{\sum_{t} c^P_t P_{0,t}}^{ \text{Real Power Cost} }
+ \overbrace{\sum_{t} c^Q_t Q_{0,t}}^{\text{ Reactive Power Cost} }
+ \overbrace{\sum_{y,t}{ c_y f_{y,t} }}^{ \text{Transformer Cost} },
\end{equation}
\emph{subject to:} network constraints \eqref{EqSubstation}--\eqref{Qinject}, transformer constraints
\eqref{Xf1}--\eqref{Cycle}, DER 
constraints \eqref{EqPVcon1}--\eqref{EVCon6}, with variables $v_{i,t}$, $l_{ij,t}$, $f_{y,t}$, $p_{s,t}$, $p_{e,t}$, $u_{e,t}$ nonnegative, and $P_{0,t}$, $Q_{0,t}$, $P_{ij,t}$, $Q_{ij,t}$, $p_{j,t}$, $q_{j,t}$, $h_{y,t}$, $q_{s,t}$, $q_{e,t}$ unrestricted in sign.

It is important to note that the quadratic equality constraint \eqref{EqCurrentDef} is non-convex. 
Following \cite{FarivarLow2013}, we relax \eqref{EqCurrentDef} to 
$v_{i,t} l_{ij,t} \geq P_{ij,t}^2 + Q_{ij,t}^2, \, \forall j,t$, and
the resulting relaxed AC OPF problem is a convex SOCP problem, which can be solved efficiently using commercially available solvers.

\section{DLMC Components} \label{DLMCs}

In this section, we provide a rigorous analysis and interpretation of DLMCs, 
which represent the dynamic marginal cost for delivering real and reactive power, denoted by P-DLMC and Q-DLMC, respectively, at a specific location and time period.
In the context of the Full-opt problem, DLMCs are obtained by the dual variables of constraints \eqref{EqRealBalance} and \eqref{EqReactiveBalance} ---
at the root node, we have $\lambda_{0,t}^P = c_t^P$, and $\lambda_{0,t}^Q = c_t^Q$.

Notably, DLMCs can be defined for any DER schedule reflected in net demand variables $p_{j,t}$ and $q_{j,t}$.
Given these variables, the AC OPF problem solution defines an operating point for the distribution network, which is obtained by the solution of the power flow equations \eqref{EqRealBalance}--\eqref{EqCurrentDef} for a given root node voltage, and is described by the real/reactive power flows, voltages and currents, i.e., variables $P_{ij,t}$, $Q_{ij,t}$, $v_{j,t}$, and $l_{ij,t}$.
DLMC decomposition employs sensitivity analysis, duality and optimality conditions of the enhanced AC OPF problem.

Let us consider DLMCs at node $j'$, time period $t'$, $\lambda_{j',t'}^P$ and $\lambda_{j',t'}^Q$.
The sensitivity of the power flow solution \emph{w.r.t.} net demand for real and reactive power, is reflected in partial derivatives $\frac{\partial P_{ij,t'}}{\partial p_{j',t'}} $, $\frac{\partial Q_{ij,t'}}{\partial p_{j',t'}} $, $\frac{\partial v_{j,t'}}{\partial p_{j',t'}}$, $\frac{\partial l_{ij,t'}}{\partial p_{j',t'}}$, and $\frac{\partial P_{ij,t'}}{\partial q_{j',t'}} $, $\frac{\partial Q_{ij,t'}}{\partial q_{j',t'}} $, $\frac{\partial v_{j,t'}}{\partial q_{j',t'}}$, $\frac{\partial l_{ij,t'}}{\partial q_{j',t'}}$, respectively. 
The calculation of sensitivities \emph{w.r.t.} real net demand at node $j'$, time period $t'$ ($4N$ unknowns) involves the solution of a linear system derived from power flow equations \eqref{EqRealBalance}-\eqref{EqCurrentDef}, after taking partial derivatives \emph{w.r.t.} $p_{j',t'}$, at time period $t'$, for all nodes ($4N$ equations).
Similarly for the sensitivities \emph{w.r.t.} reactive net demand.
In total, it requires the (amenable to parallelization) solution of $2NT$ linear systems (referring to sensitivities \emph{w.r.t.} real and reactive power net demand, calculated at $N$ nodes, for $T$ time periods), each system involving $4N$ equations and $4N$ unknowns.

The DLMC decomposition to additive components follows:
$\lambda_{j',t'}^P 
= c_{t'}^P 
+ \overbrace{c_{t'}^P  \sum_j r_{ij}\frac{\partial l_{ij,t'}}{\partial p_{j',t'}}}^{\text{Real Power Marginal Losses}}
+ \overbrace{c_{t'}^Q  \sum_j x_{ij}\frac{\partial l_{ij,t'}}{\partial p_{j',t'}}}^{\text{Reactive Power Marginal Losses}}$ \\
$ + \overbrace{\sum_j{ \mu_{j,t'} \frac{\partial v_{j,t'}}{\partial p_{j',t'}}}}^{\text{Voltage Congestion}} 
+ \overbrace{\sum_j { \nu_{j,t'} \frac{\partial l_{ij,t'}}{\partial p_{j',t'}}}}^{\text{Ampacity Congestion}}
+ \overbrace{\sum_{y}{\pi_{y,t'} \frac{\partial l_{y,t'}}{\partial p_{j',t'}}}}^{\text{Transformer Degradation}}$,\\
$\lambda_{j',t'}^Q  
= c_{t'}^Q 
+ \overbrace{c_{t'}^P  \sum_j r_{ij}\frac{\partial l_{ij,t'}}{\partial q_{j',t'}}}^{\text{Real Power Marginal Losses}}
+ \overbrace{c_{t'}^Q  \sum_j x_{ij}\frac{\partial l_{ij,t'}}{\partial q_{j',t'}}}^{\text{Reactive Power Marginal Losses}}$ \\
$ + \overbrace{\sum_j{ \mu_{j,t'} \frac{\partial v_{j,t'}}{\partial q_{j',t'}}}}^{\text{Voltage Congestion}} 
+ \overbrace{\sum_j { \nu_{j,t'} \frac{\partial l_{ij,t'}}{\partial q_{j',t'}}}}^{\text{Ampacity Congestion}}
+ \overbrace{\sum_{y}{\pi_{y,t'} \frac{\partial l_{y,t'}}{\partial q_{j',t'}}}}^{\text{Transformer Degradation}}$,
where parameter $\pi_{y,t'}$ 
includes the intertemporal impact of the transformer degradation component.
For both expressions, the first three terms are obtained by associating equations \eqref{EqRealBalance} and \eqref{EqReactiveBalance} recursively to \eqref{EqSubstation}, and taking the partial derivatives.
The first component is the marginal cost of real/reactive power at the root node.
The second and the third components represent the contribution of real and reactive power marginal losses that are sensitive to changes in net real and reactive power demand.
The fourth component reflects voltage congestion, with $\mu_{j,t'} = \bar \mu_{j,t'} - \underline \mu_{j,t'}$.
It is nonzero only at nodes with a binding voltage constraint, i.e., when $\bar \mu_{j,t'} > 0$ or $\underline \mu_{j,t'} > 0$. 
The fifth component represents ampacity congestion.
It is non zero only at lines with a binding ampacity constraint, i.e., when $\nu_{j,t'} > 0$.
Voltage and ampacity congestion components are derived by appending the active (binding) constraints in the Full-opt Lagrangian;  
alternatively, they can be obtained by the optimality conditions. 
Last, and perhaps most important, the sixth component represents the impact on transformer degradation costs that are intertemporally coupled with real and reactive power injections during preceding time periods. 
We next derive the formula for $\pi_{y,t'}$ and elaborate on it.

Consider, without loss of generality, the P-DLMC component.
Applying \eqref{Xf2} recursively, we get 
$h_{y,t} = \delta^t h_{y,0} + \sum_{\tau = 1}^t \delta^{t-\tau} \left( \epsilon_y l_{y,\tau} + \zeta_{y,\tau} \right)$, $\forall y, t$.
Taking the partial derivatives of $h_{y,t}$ \emph{w.r.t.} $p_{j',t'}$, we get:
$\frac{\partial h_{y,t}}{\partial p_{j',t'}} = \epsilon_y \delta^{t-t'}  \frac{\partial l_{y,t'}}{\partial p_{j',t'}} $, $\forall y, t \geq t'$.
From the binding transformer constraints \eqref{Xf1}, using optimality conditions and $\frac{\partial h_{y,t}}{\partial p_{j',t'}}$, 
the impact within the optimization horizon is given by:
\begin{equation*}
\begin{split}
\sum_{y,t}{c_y \frac{\partial f_{y,t}}{\partial p_{j',t'}}} 
= &\sum_{y,t, \kappa}{\xi_{y,t,\kappa} \big( \alpha_{\kappa} \frac{\partial h_{y,t}}{\partial p_{j',t'}} + \beta_{y, \kappa} \frac{\partial l_{y,t}}{\partial p_{j',t'}}  \big)}
= \\
&\big( \epsilon_y \sum_{t=t'}^T { \delta^{t-t'} \tilde \xi_{y,t}  } 
+ \sum_{\kappa} \xi_{y,t',\kappa} \beta_{y,\kappa} \big) \frac{\partial l_{y,t}}{\partial p_{j',t'}},
\end{split}
\end{equation*}
where $\tilde \xi_{y,t} = \sum_{\kappa} \xi_{y,t,\kappa} a_{\kappa}$.
We note that $\sum_{\kappa}\xi_{y,t,\kappa} = c_y$, and, obviously, the sum  $\sum_{\kappa} \xi_{y,t,\kappa}$ involves only the binding constraints (otherwise $\xi_{y,t,\kappa} = 0$); if only one constraint is active, then the respective dual $\xi_{y,t,\kappa}$ should equal $c_y$.
Adding the impact that extends beyond the optimization horizon, which is obtained by appending \eqref{Cycle} in the Lagrangian yielding $\rho_y \frac{\partial h_{y,T}}{\partial p_{j',t'}} = \epsilon_y \delta^{T-t'}  \rho_y \frac{\partial l_{y,t'}}{\partial p_{j',t'}}$,
$\pi_{y,t'}$ is given by:
\begin{equation} \label{piyt}
\pi_{y,t'} =  \epsilon_y \Big( \sum_{t=t'}^T  \delta^{t-t'} \tilde \xi_{y,t} +  \delta^{T-t'}  \rho_y \Big)
+ \eta_y \tilde \xi_{y,t'},
\end{equation}
where  
$\eta_y = \frac{4 \Delta \bar{\theta}_{y}^H}{ 5 l_y^N}$
(hence, $\beta_{y,\kappa} = \eta_y \alpha_{\kappa}$, and $\sum_{\kappa} \xi_{y,t',\kappa} \beta_{y,\kappa} = \eta_y \tilde \xi_{y,t'}$).
Notably, $\pi_{y,t'}$ involves a summation over time that is related to the transformer thermal response dynamics and captures subsequent time period costs.
The first term in \eqref{piyt} that is multiplied by $\epsilon_y$ refers to the contribution of the load at time period $t'$ to the TO temperature at $t'$ and subsequent time periods as well as the impact that extends beyond the optimization horizon (involving dual variable $\rho_y$), whereas the second term that is multiplied by $\eta_y$ refers to the contribution to the winding temperature rise over the TO.
Recall that the winding time constant is small compared to our time period length (1 hour), hence its impact is limited to time period $t'$, whereas the oil time constant is larger, hence the intertemporal impact on subsequent time periods is not negligible.
This intertemporal impact is discounted by $\delta^{t-t'}$ (with $\delta < 1$), but it also considers the slope of transformer LoL $a_{\kappa}$ and applies a higher weight to high HST time periods --- recall that this slope is large when the HST is high reflecting current and past period transformer loading.
The impact that extends beyond the optimization horizon is discounted by $\delta^{T-t'}$, hence, it may become significant for hours closer to the end of the horizon.

Lastly, we note that the decomposition of DLMCs to additive components enables the rigorous quantification of costs by source and more importantly the use of distributed algorithms to allow the scaling of AC OPF to real distribution feeders.
The contribution of each component can primarily provide useful information for pricing purposes, the potential allocation of associated revenue and for network asset and DER capacity expansion when the marginal costs obtained at the optimal DER schedule are sufficiently high. 
In addition, however, the fact that we can obtain marginal-cost-component-specific sensitivities for a given system operating point allow us, in the absence of active congestion constraints, to estimate the marginal costs without solving the AC OPF, but using instead a power flow solver.
We note that the transformer marginal cost component can be estimated by calculating the temperature trajectory and using the aging model.

\section{Case Study Based Numerical Findings} \label{CaseStudy}

In this section, we present the case study input data and scheduling options (Subsection \ref{FeederData}),
we discuss the main findings (Subsection \ref{NumResults}),
and we elaborate on the transformer degradation intertemporal impact (Subsection \ref{intertemp}).

\subsection{Case Study Input Data and Scheduling Options} \label{FeederData}

\begin{figure}[b]
\centering
\includegraphics[width=3.45in]{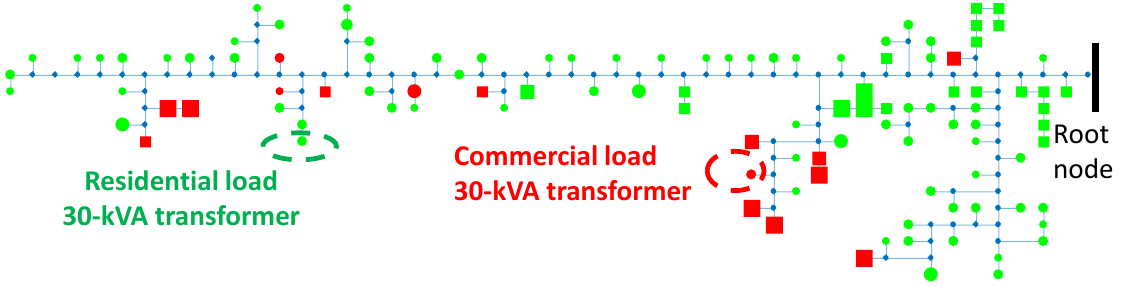}
\caption{Feeder topology diagram. Transformers: round (pole); square (pad); green (residential); red (commercial); sizes reflect nameplate capacity. Feeder data (lines, transformers, load profiles, irradation level, ambient temperature) are listed in \cite[Section III.A]{PartII}.}
\label{figFeeder}
\end{figure}
Our proposed model is illustrated on a 13.8-kV, 307-node, 110-transformer, feeder of Holyoke Gas and Electric (HGE), a municipal distribution utility in MA, US. 
To improve the exposition of results while emphasizing the local effect of DERs, we consider different levels of EVs and PVs connected exclusively to two nodes, depicted in Fig. \ref{figFeeder}, representing two 30-kVA transformers serving commercial and residential loads.
For both transformers, $R=5$, $\Delta \bar \theta^{TO} = 55$, 
$\Delta \bar \theta^{H} = 25$, and $M=8$ ($\tilde{F}_{AA}$ linearization follows \cite[Fig. 1]{PartI};
indicatively, with 3 equal segments from $110^o$C to $140^o$C the maximum error is less than 1 hour, whereas 6 equal segments would reduce the maximum error by a factor of 3 overestimating LoL by at most 3\%.
At the commercial (residential) node, EVs are connected 9am--5pm (7pm--7am) and require charging 12 kWh (18 kWh).
At departure, EVs must be fully charged.
EV battery capacity is 50 kWh, the maximum charging rate 3.3 kW/h, and the charger capacity 6.6 kVA.
PVs are 10-kVA rooftop solar units.
To model initial conditions in the daily cycle reasonably, transformer temperatures were required to coincide at the beginning and at the end of the cycle; a similar constraint was imposed on EV battery SoC.
\begin{table*}[ht] 
	\caption{Aggregate System Cost Difference (in \$) and LoL (in hours)}  \label{tab3a} 
\centering
\includegraphics[width=6.5in]{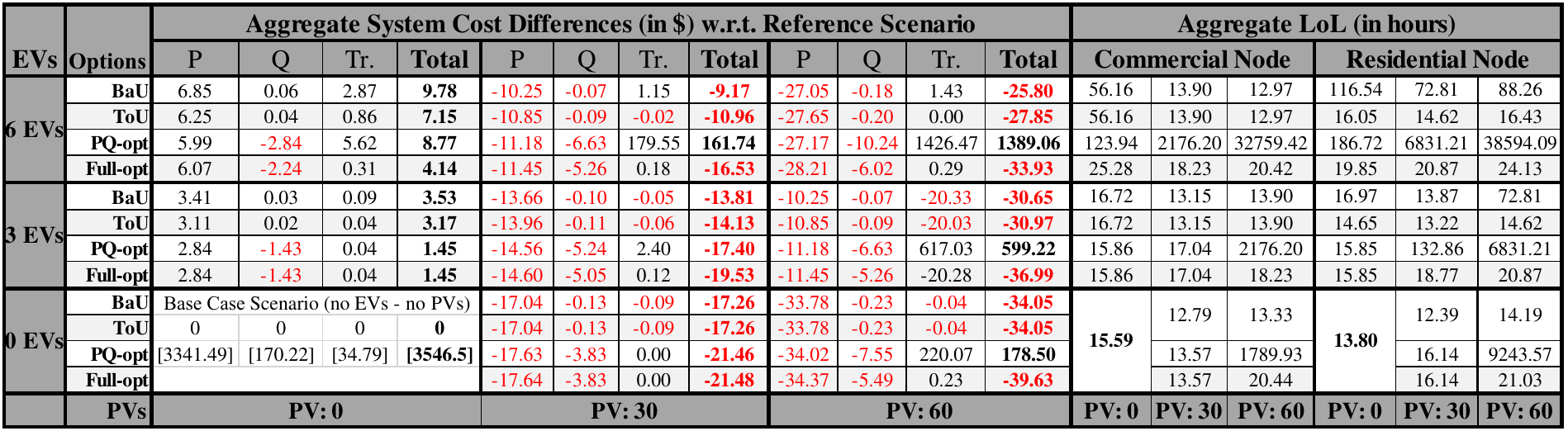}
\end{table*} 

We consider the following 4 EV/PV scheduling options:\\
$\bullet$ ``\textbf{BaU}'' (Business as Usual): EVs ``dumb'' charge at full rate upon arrival with p.f. = 1. 
PVs operate with p.f. = 1.\\
$\bullet$ ``\textbf{ToU}'' (Time-of-Use): EVs charge to minimize real power cost relative to the selected LMP trajectory while maintaining p.f. = 1. PVs operate with p.f. = 1.\\
$\bullet$ ``\textbf{PQ-opt}'': EVs/PVs are scheduled in order to minimize real and reactive power cost, subject to voltage and ampacity constraints while ignoring transformer degradation cost.
PQ-opt can be thought of as a traditional line loss minimization problem akin to the day-ahead problem presented in \cite{BaiEtAl-DLMP}.
It can be alternatively thought of as a special case of Full-opt below without considering transformer degradation costs.\\
$\bullet$ ``\textbf{Full-opt}'': EVs/PVs are scheduled by adapting to full DLMCs including transformer degradation costs, given by the solution of Full-opt summarized in \eqref{Obj2}.\\

\subsection{Numerical Results} \label{NumResults}
Optimization problems were solved on a Dell Intel Core i7-5500U @2.4 GHz with 8 GB RAM, using CPLEX 12.7, solution times were up to 10 sec, and convex relaxations were exact in all instances.

In Table \ref{tab3a}, we compare several DER penetration scenarios at the nodes of interest:
EVs equal 0, 3, and 6, and rooftop PV installations equal 0, 30, and 60-kVA (i.e., 0, 3, and 6 units of 10-kVA rooftop solar, respectively).
The scenario with 0 EVs and PVs serves as the base case for reference purposes.
Table \ref{tab3a} shows the system-wide daily cost differences (in \$) for real/reactive power (P/Q), transformer degradation, and total cost relative to the base case.
For readability, negative values are shown in red.
Table \ref{tab3a} reports also the aggregate daily Loss of Life (LoL) of the two 30-kVA transformers.
As expected, Full-opt achieves the lowest total cost, since it co-optimizes real/reactive power and transformer cost and maintains low LoL.
PQ-opt achieves the lowest combined real and reactive power costs (P and Q columns in Table \ref{tab3a}), but exhibits some very high values of aggregate LoL and transformer degradation cost, as a result of high reactive power provision.
Unlike PQ-opt and Full-opt, ToU and BaU do not take advantage of DER reactive power provision capabilities, but in some cases, particularly in the residential node, achieve low aggregate LoL.
 
In the commercial node, BaU and ToU produce identical EV schedules (by coincidence LMPs are increasing during the day).
In the residential node, ToU exhibits a 3-hour shift relative to the BaU EV profile (EVs start charging at 10pm). 
For the specific LMPs and EV/PV penetrations, the ToU schedule results in sustainable transformer LoL performance;
different LMPs though and/or higher penetration levels may result in worse performance of the ToU option. 
Since we focus on understanding the granularity of DLMCs, we did not consider multiple EV profiles and arrival times that would demonstrate the potential synchronization under the ToU schedule (see e.g., the simulations in \cite{HilsheyEtAl_2013}) resulting in worse performance relative to the BaU schedule.
For the high EV/PV penetration presented later, however, high LoL is associated with ToU.

It is important to note that with the exception of Full-opt generated schedules, all other DER schedules are not adapted to DLMCs.
Under BaU, ToU, and PQ-opt, DLMCs simply represent the marginal cost of delivering an incremental unit of real and reactive power to/at a specific network location and time, conditional upon the given DER (EV/PV) schedules. 
As such, DLMCs under non-adaptive schedules should be interpreted as ex post spatiotemporal marginal costs indicating desirable directions of change.
So for a given DER schedule (BU, ToU, PQ-opt), DLMCs are calculated solving a partial Full-opt problem with given/fixed EV/PV schedules.
As such, DLMCs can be interpreted as economically/socially efficient financial incentives indicating a direction for improving system costs.
For example, a DLMC at a specific node and time, which is higher than the Full-opt DLMC, signals that less power should be consumed then and there, if at all possible.

\begin{figure}[tb]
\centering
\includegraphics[width=3.45in]{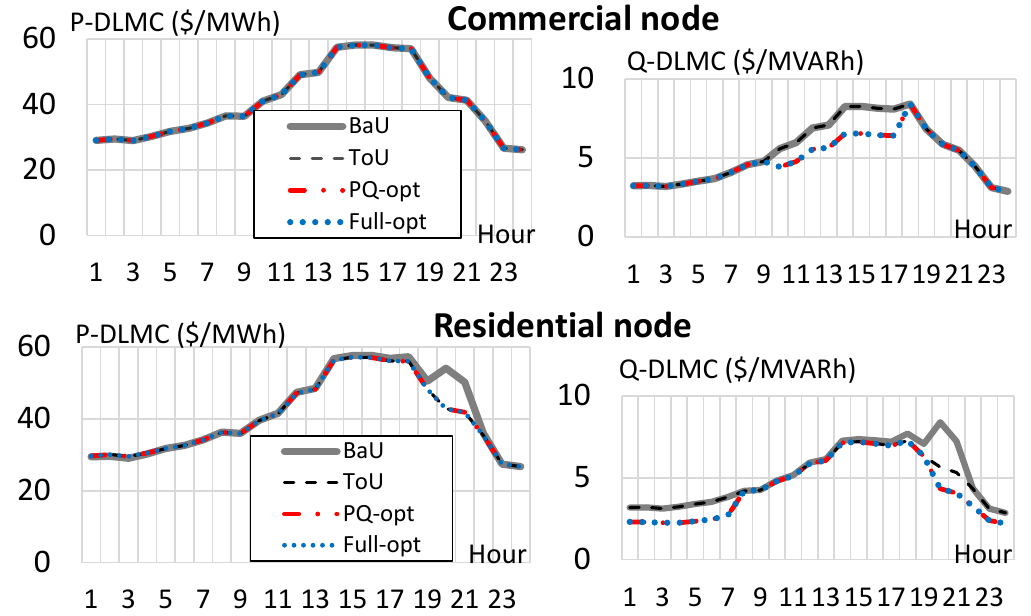}
\caption{DLMCs, commercial and residential nodes, 3 EVs.}
\label{figEV3A}
\end{figure}
Let us first consider the EV-only penetration scenarios.
In Fig. \ref{figEV3A}, we show P-DLMCs and Q-DLMCs for a relatively low 3-EV scenario. 
The differences are rather small and are related to difference in the provision of reactive power.
At the commercial node, P-DLMCs are practically the same under all scheduling options. 
Q-DLMCs are higher under BaU/ToU compared to PQ-opt/Full-opt, indicating that higher reactive power provision would be desirable.
At the residential node, P-DLMCs associated with BaU EV schedules are higher for hours 20--22, pointing to the desirability of a 3-hour shift in EV consumption;
this coincides with the result of other scheduling options (EVs start charging at 10pm).
Q-DLMCs associated with BaU/ToU EV schedules suggest the desirability of providing reactive power;
this coincides with the results of Full-opt and PQ-opt that indeed provide reactive power and decrease system-wide cost.
P-DLMC components range from 1.7\% to 7.5\% for real power losses, from 0.2\% to 0.9\% for reactive power losses and from 0.2\% to 1.7\% for transformer degradation. 
For Q-DLMCs, during hours 10--17 (when EVs are plugged in), the marginal cost components for real and reactive power are lower under Full-opt, while the transformer degradation component becomes negative.
This is interestingly associated with a reverse reactive power flow at that node and those time periods rendering the sensitivity of the current \emph{w.r.t.} reactive power net demand negative. 
Q-DLMC components under Full-opt range from 8\% to 30.7\% for real power losses, from 1\% to 3.7\% for reactive power losses, and from -2\% to 5\% for transformer degradation.

\begin{figure}[tb]
\centering
\includegraphics[width=3.45in]{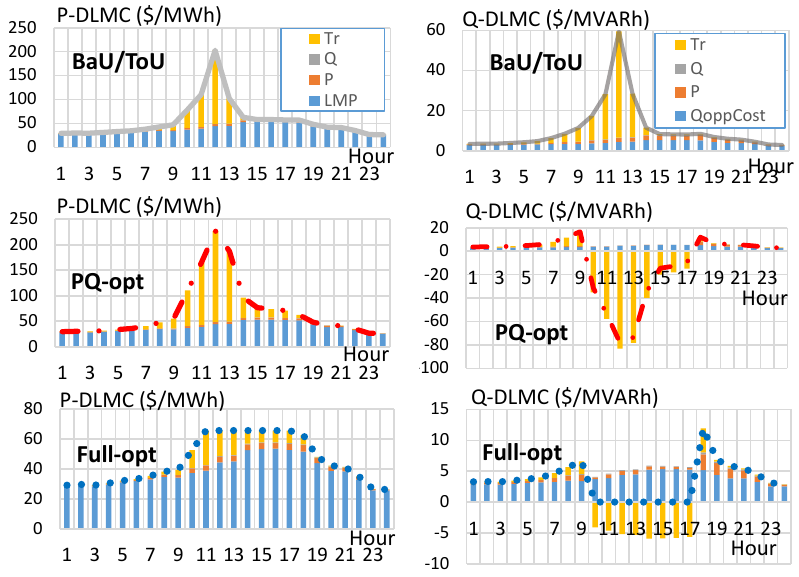}
\caption{DLMC components, commercial node, 6 EVs.}
\label{figEV6D}
\end{figure}
The 6-EV penetration scenario is more impactful, and PQ-opt is associated with a 5 to 9 fold higher aggregate transformer LoL relative to Full-opt (see Table \ref{tab3a}).
For this scenario, Fig. \ref{figEV6D} shows DLMCs amd their components at the commercial node.
P-DLMCs exhibit spikes under BaU/ToU and PQ-opt implying that charging rates are too high, whereas Full-opt exhibits a smoother profile.
Q-DLMCs exhibit positive spikes under BaU/ToU, implying that EVs should provide reactive power. 
Negative spikes under PQ-opt suggest that EVs provide excessive reactive power. 
We note that Full-opt results in practically zero Q-DLMCs ($\to 0^+$) while EVs are plugged in, thus implying high reactive power availability. 
Nevertheless, despite the fact that the imputed reactive power income --- 
assuming self-scheduling EVs considering DLMCs as prices --- 
decreases with tanking Q-DLMCs, Full-opt results in lower P-DLMCs as well; the net impact on the overall system cost from EV charging is lower and EVs realize less imputed income from Reactive Power ``sales'' but overall their effective charging costs are lower!
Full-opt Q-DLMCs support the system-optimal solution by incentivizing EVs to provide reactive power at a rate that is lower than their charger capability,
whereas PQ-opt ``brute-forces'' EVs to fully utilize their inverter capability resulting in negative Q-DLMCs that significantly increase the imputed cost of EVs. 

It is noteworthy that the marginal transformer degradation cost varies across scenarios and scheduling options illustrating the DLMC-implied incentives.
Indeed, the DLMC spikes (positive/negative) are caused by the transformer degradation component.
In general, low transformer costs incurred by Full-opt across all hours are associated with smoother DLMC profiles.
Notably, although voltage and ampacity DLMC components are zero, i.e., the respective constraints are not binding, transformer degradation is high. 
This is an interesting observation showing that traditional line loss minimization approaches employing Volt/VAR control capabilities of smart inverters and capacitors, such as in PQ-opt, may ``mask'' the significant impact of transformer degradation on DLMCs, and impose potentially excessive cost on the distribution network assets.
Similar remarks can be made for the residential node. 

\begin{figure}[tb]
\centering
\includegraphics[width=3.45in]{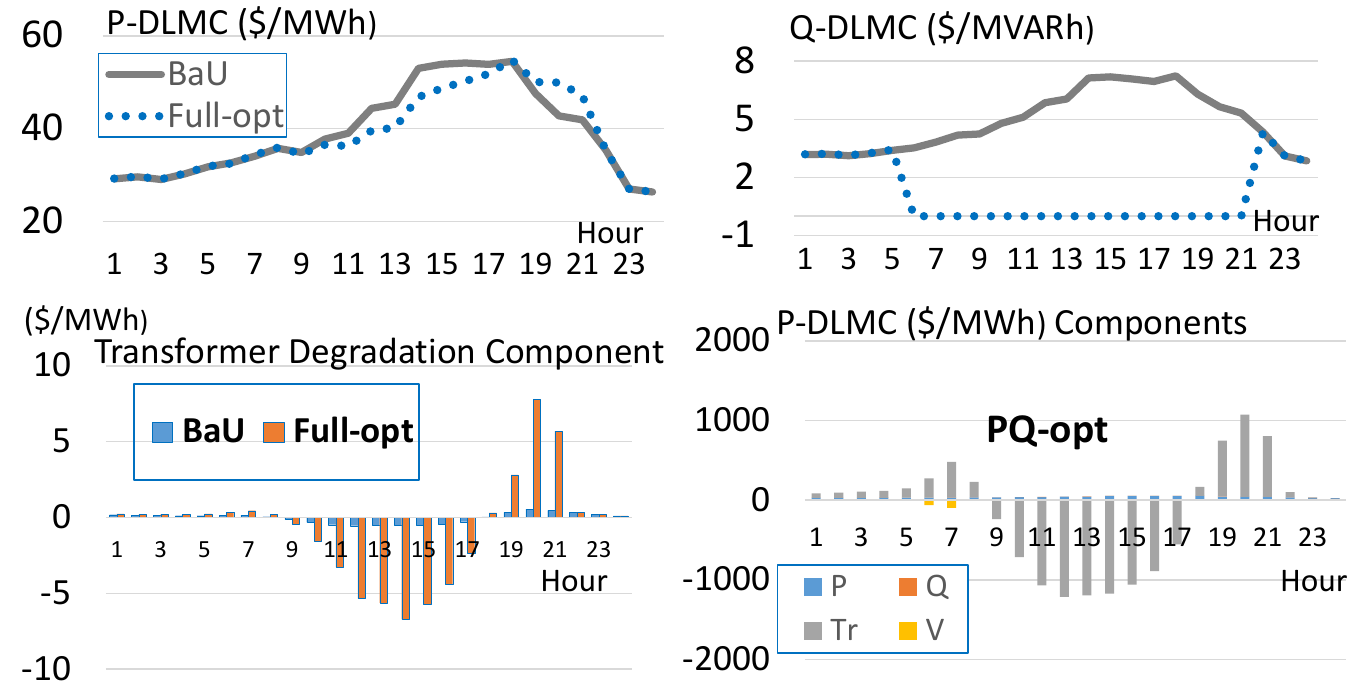}
\caption{DLMCs, residential node, 60-kVA PV.}
\label{figPV60}
\end{figure}
Considering the PV-only penetration scenarios, Table \ref{tab3a} suggests that the 30-kVA PV scenario has rather insignificant impact on network costs, whereas the 60-kVA scenario is quite impactful, with PQ-opt yielding very high aggregate LoL values compared to BaU or Full-opt. 
For the latter scenario, at the residential node, we illustrate, in Fig. \ref{figPV60}, P-DLMCs and Q-DLMCs under BaU and Full-opt (top figure), P-DLMC transformer degradation components for BaU and Full-opt (bottom left), and P-DLMC components for PQ-opt (bottom right).
We observe that BaU and Full-opt result in similar P-DLMCs (top left); the differences are explained by the transformer degradation component (bottom left).
For BaU and Full-opt, the transformer component becomes negative during hours 9--17, when there is a reverse real power flow (PV generation exceeds the load) and the sensitivity of the current \emph{w.r.t.} net demand becomes negative.
The transformer component under PQ-opt (bottom right) is ``amplified'' --- at least two orders of magnitude higher than Full-opt.
We also observe a negative voltage component under PQ-opt at hours 6 and 7, when the voltage upper bound is reached, indicating a negative sensitivity of the voltage \emph{w.r.t.} net demand.
Indeed, increasing net demand would result in a voltage drop.
However, the transformer degradation component at those hours is positive, indicating a positive sensitivity of the current \emph{w.r.t.} net demand, and is about 4 to 5 times larger in magnitude than the voltage component.
Similar to the 6-EV scenario, DLMC differences are caused by reactive power provision.
Positive Q-DLMCs under BaU (top right) encourage the provision of reactive power.
Under PQ-opt, Q-DLMCs (not shown in Fig. \ref{figPV60}) become highly negative, suggesting that reactive power provision is excessive (PVs fully utilize their inverter capabilities).
Under Full-opt, PVs provide reactive power as needed to minimize system cost; 
Q-DLMCs tank, indicating that there is excess installed capacity and it is optimal to underutilize it.
Full utilization, as is the case under PQ-opt, renders negative Q-DLMCs and interacts with P-DLMCs to result in a large transformer degradation component.

Considering the EV-PV synergy, with the exception of PQ-opt at the residential node, the 3-EV and 30-kVA PV scenarios are easily sustainable with low aggregate LoL. 
At hours 6, 7, 20, and 21, high reverse reactive power flows overload the transformer. 
P-DLMCs exhibit positive spikes and Q-DLMCs negative spikes at those hours. 
Increasing EVs to 6 and/or PVs to 60-kVA results in very high LoL under PQ-opt.
On the other hand, as expected, Full-opt achieves sustainability with low LoL under all scenarios.

Lastly, we consider increased DER penetrations.
Full-opt is still able to accommodate 9 and 12 EVs, whereas other scheduling options fail.
ToU keeps reasonably low aggregate LoL (40 hours) only for the 9-EV scenario at the residential node, although still 2-fold higher than Full-opt, whereas 32-fold higher LoL values are observed for the 12-EV scenario.
At the commercial node, BaU/ToU exhibit aggregate LoL values that exceed 1000 hours (PQ-opt is even worse).
The robust performance of Full-opt with increasing EV penetrations is explained by an appropriate reduction in the charging rate and the utilization of the inverter. 
Notably, for the 12-EV scenario, Full-opt reaches 114 aggregate LoL hours, which is still 40 times lower than the values encountered under BaU/ToU schedules.
Experimentation with 90 and 120-kVA PV penetration scenarios showed that Full-opt can significantly increase PV hosting capacity ``curtailing'' real and reactive power as needed, and, thus, avoiding excessive transformer temperatures and LoL.
Notably, if we consider self-scheduling PVs responding to optimal DLMC-based signals, Full-opt DLMCs support economically efficient ``curtailment.'' 

\subsection{Transformer Degradation Intertemporal Impact} \label{intertemp}

In general, the transformer degradation component at a specific location includes the impact on all transformers through weights that are proportional to the sensitivities of the current \emph{w.r.t.} net demand.
Since in our feeder, overloads are located on two rather distant transformers, we focus our attention on the impact only on the transformer located at a specific node.
We also note that a numerical evaluation of the thermal model approximations in the HST calculations yielded very accurate results with average errors of about 1\%-2\% or lower in the high loading conditions of interest to our analysis.

Let us consider transformer $y$ denoted by line $ij'$. 
Its P-DLMC degradation component at node $j'$, time period $t'$, $\pi_{y,t'}  \frac{\partial l_{y,t'}}{\partial p_{j',t'}}$, is given by
$ 
\big[\big(\frac{\partial l_{y,t'}}{\partial p_{j',t'}}\big) / l_y^N \big]
\big\{ \rho_y \big[ \frac{55}{6} { \big( \frac{3}{4} \big)}^{T-t'}  
 \big]  
 +    c_y  
 \big[
 \frac{55}{6}  \tilde \alpha_{y,t'} 
 + 20 \tilde \alpha_{y,t'}
 + \sum_{\tau = 1}^{T-t'} { \left( \frac{3}{4} \right)}^{\tau} \tilde \alpha_{y,t'+ \tau}
 \big] \big\}
 $,
where $\rho_y$ captures the impact of increasing HST in hour 24 on future HST and transformer degradation,
$\tilde \alpha_{y,t} = \frac{1}{c_y} \tilde \xi_{y,t} = \frac{1}{c_y}\sum_{\kappa} \xi_{y,t,\kappa} a_{\kappa} $ is an adjusted slope of the aging acceleration factor.
First, we note that the component depends on the sensitivity of the current \emph{w.r.t.} net demand at hour $t'$, which determines the sign of the component, normalized by the transformer nominal current.
Second, the component depends on the slope of the aging acceleration factor at hour $t'$, as well as the respective slopes of subsequent hours.
The impact that is related to hour $t'$ includes the contribution of the TO temperature and the winding.
In fact, there is a fixed ratio between the two that depends on the transformer parameters, which in our case is equal to 2.18 (the contribution of the winding is higher).
The impact that is related to subsequent periods refers only to the TO temperature dynamics.
Interestingly, although this impact decays with ${ \left( \frac{3}{4} \right)}^{\tau}$, for each hour $\tau$ that follows $t'$, it can become significant when the slopes $\tilde \alpha_{y,t'+ \tau}$ are quite large.
As we show next, the impact that extends beyond the optimization horizon that is related to $\rho_y$ may also become significant.

\begin{figure}[tb]
\centering
\includegraphics[width=3.4in]{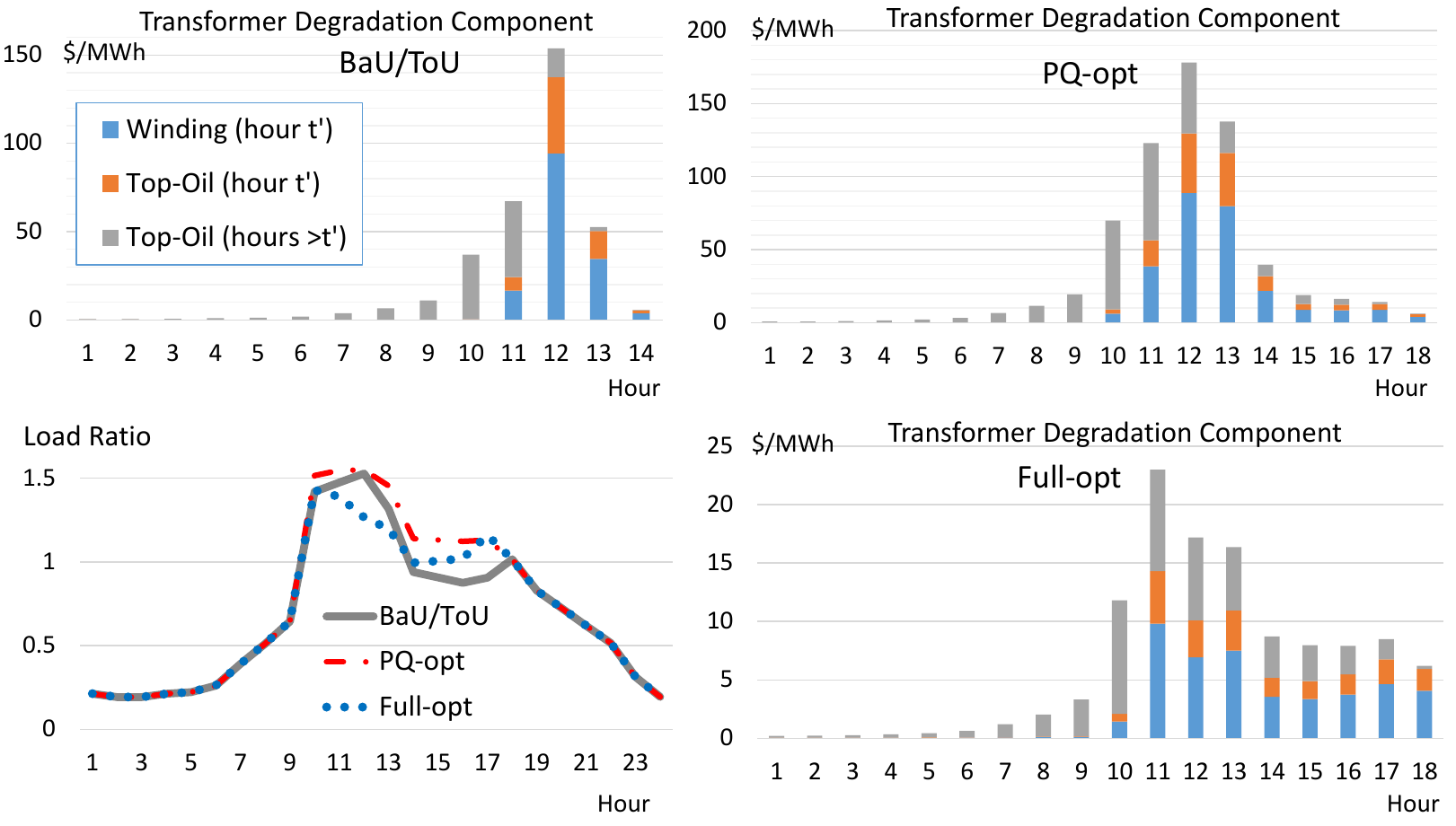}
\caption{P-DLMC transformer degradation components, commercial node, 6 EVs. Load Ratio $l_{y,t}/l_y^N$.}
\label{figEV6E}
\end{figure}
Let us consider the P-DLMCs for the 6-EV penetration scenario, at the commercial node, shown in Fig. \ref{figEV6D} (left).
For this scenario, we show in Fig. \ref{figEV6E} the load ratio and the transformer degradation sub-components.
An interesting remark is that the transformer component becomes significant before the EVs plug in, e.g., at hour 9, due to the TO temperature impact on subsequent hours. 
Under BaU/ToU schedules, the transformer component becomes 10.96 (\$/MWh) at hour 9, accounting for 23.3\% of the P-DLMC.
The component is still significant even at hours 7 and 8, accounting for 9.95\% and 15.3\% of the P-DLMC, respectively.
It becomes as high as 153.75 (\$/MWh) at hour 9, and diminishes after hour 14 as the transformer cools down.
At hours that precede the overloading hours (due to EVs) and at the first two hours of charging, the impact on the subsequent hours dominates the transformer component.
The winding and TO impact become dominant at hour 12.
PQ-opt exhibits similar results which are amplified and extended (in time duration) by the excessive brute-force-driven reactive power provision.
The transformer component becomes significant starting at hour 6 (9.3\% of the P-DLMC), it reaches a maximum of 178.1 (\$/MWh) at hour 12, the impact of subsequent hours remains large until hour 11, and then diminishes after hour 18.
Similar intertemporal effects are observed under the Full-opt, although the behavior is much smoother compared to PQ-opt.
The impact on subsequent hours makes up the major part of the P-DLMC until hour 10, and even though the transformer component becomes higher relative to BaU/ToU at hours 14--18 (ranging from 3 to 8 \$/MWh), the profile is smoother.

\section{Conclusions and Further Research} \label{Conclusions}

We have provided an enhanced AC OPF for granular marginal costing on distribution networks including asset degradation costs exemplified by transformer marginal LoL.
We have also derived and analyzed additive DLMC components, and illustrated the intertemporal characteristics of the transformer degradation component.
The DLMCs estimated by the proposed AC OPF are compatible with optimal DER self-scheduling that is fully adapted to the optimal distribution network configuration. 
As such, the DLMCs represent the spatiotemporal marginal costs on the distribution network that co-optimize network and DER scheduling decisions. 
Whereas we have considered the real and reactive power marginal cost at the substation bus to be given, 
it is possible that after all distribution networks and their DER schedules are co-optimized, the Wholesale Transmission market adjusts producing new LMPs.
Although LMPs are most likely to be insensitive to real power DLMCs that affect primarily distribution node injections but less so the aggregate T\&D power exchange, it is possible that reactive power compensation at the T\&D interface is impacted by the distribution network reactive power schedule.
This is a topic of future work that may require ISOs to provide reasonable, possibly extended LMP type, price functions at the T\&D interface.

Our numerical results inquiry has focused on service transformers that are sensitive to granular/location-specific DER schedules.
A study focusing on the T\&D interface particularly during the few weeks following a station transformer failure, could and should model station transformers placing them on the line connecting the root node to the immediately downstream node.
During periods of station transformer overloading, the root node LMP would be translated to a higher downstream DLMC.

Employing an actual distribution feeder, we have shown that optimal DLMC-driven DER scheduling can achieve significantly higher hosting of DERs by distribution networks with minimal infrastructure investments.
In a practical application that one may expect to be adopted in the future, DERs that opt in a DLMC-based rate would most likely see a two-part tariff, consisting of a DLMC part and a constant part aimed to address fixed cost recovery issues.
We anticipate that sophisticated rate designs that fine tune two-part tariff ideas will attract the attention of capable system engineers, economists and practitioners.
Most importantly, we have shown that when DER schedules ``adapt'' to feeder DLMCs, the equilibrium DLMC profiles are rather smooth, thus allaying fears that DLMC volatility may lead to ratepayer revolt if DLMC-based rates are adopted.
Comparison with popular open-loop DER scheduling options provided solid evidence that optimal DLMC-based clearing markets can bring about significant economic efficiencies and support the sustainability/adequacy of current distribution network infrastructure in the presence of high DER (PV, EV, and the like) adoption. 
It is also particularly notable that the presence of DERs scheduled optimally through spatiotemporal DLMC adaptation will have a significant impact on grid asset capacity expansion choices. 

In this work, the optimal daily/24-hour operational planning task employed a centralized decision maker who is aware of both network information and DER capabilities, preferences and constraints.
Alternatively, decentralized, distributed algorithms \cite{MolzahnEtAl_2017} can be employed, in which DERs respond to hour and location specific DLMCs communicated to them by the network operator to derive tentative schedules and interact iteratively until the optimal daily operational schedules are obtained.
In addition, we note that a distributed algorithm approach to operational scheduling would be preferable since it would enable individual DERs to respond to day-ahead DLMCs by optimizing the sum of operational scheduling benefits plus likely real-time preference variations, however large or small but likely to vary across individual DERs, as represented by the uncertain real-time DLMCs.

Furthermore, we should acknowledge that the role of uncertainty will become significant as more stochastic resources are coming into play. 
The well-known remedy has been to secure reserves;
in fact \cite{CaramanisEtAl2016} is the first DLMP work that discusses reserves in distribution networks. Notably, defining the reserve requirements is another open research area. 
We envision that in the future, reserves will be quite likely defined endogenously replacing rules-of-thumb and offline studies that determine mostly static and generally overly conservative security margins in current utility practice.
In \cite{MiethDvorkin_2019}, a chance-constrained AC OPF formulation is employed to introduce uncertainty in radial distribution networks, considering only small-scale generators as controllable DERs, and treating all behind-the-meter DERs as uncontrollable.
However, we envision that the latter DERs have a great potential to be controlled and scheduled, thus playing a key role in dealing with uncertainty in distribution networks, and we believe that this is a research direction that is certainly worth pursuing.

Lastly, our ongoing research includes extending our DLMC framework to a 3-phase network representation, considering unbalanced distribution grids, addressing degradation of other network assets, e.g., distribution lines, and employing decomposition approaches and distributed algorithms \cite{Allerton2019} that rely on our unbundling of DLMCs to additive components so that we can deal with real size multiple feeder systems and large numbers of diverse DERs transcending EVs and PVs and including micro generators, smart buildings with pre-cooling/heating capable HVAC, smart appliances, storage, and new technologies.

\section*{Acknowledgment}
The authors would like to thank Holyoke Gas and Electric for providing actual feeder data for the pilot study.

\ifCLASSOPTIONcaptionsoff 
  \newpage
\fi

\begin{IEEEbiography}
[{\includegraphics[width=1in,height=1.25in,clip,keepaspectratio]{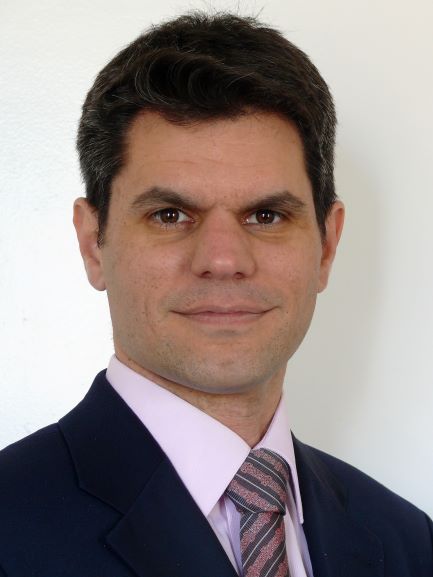}}]
{Panagiotis Andrianesis} is a graduate of the Hellenic Army Academy, also
holding a B.Sc. degree in economics (2004) from the National and Kapodistrian
University of Athens, and a Diploma degree in electrical and computer
engineering (2010) from the National Technical University of Athens, Greece.
He received his M.Sc. degree in production management (2011) and his Ph.D. degree
(2016) in the area of design and analysis of electricity market mechanisms from the
University of Thessaly (UTh), Greece.
He is currently a Postdoctoral Associate in the Division of Systems Engineering, at
Boston University, also affiliated with the Information and Data Science Research
Group. He has been a Research Associate of the Production Management Laboratory, at the UTh,
and a Consultant and Research Associate of ECCO International Inc. 
His research interests include power system economics, electricity markets,
operations research, optimization, and applied mathematics.
Dr. Andrianesis is a Member of IEEE, INFORMS, MOS, and POMS. He is the recipient of
the 2010 IEEE APS Pre-Doctoral Research Award.
\end{IEEEbiography}

\begin{IEEEbiography}
[{\includegraphics[width=1in,height=1.25in,clip,keepaspectratio]{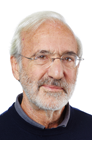}}]
{Michael Caramanis} received the B.S. in Chemical Eng. from Stanford University, Palo
Alto, CA, USA, in 1971 and the M.S. and Ph.D. degrees in Engineering from Harvard
University, Cambridge, MA, USA, in 1972 and 1976, respectively. Since 1982 he is at
Boston University where he is Professor of Systems and Mechanical Eng. He chaired the
Greek Regulatory Authority for Energy and the International Energy Charter's
Investment Group (2014-2008), was personally involved in power market implementations
in England (1989-90) and Italy (2000-03), and his written work has influenced Power
Market design in the U.S. and Europe. His current application domain focus is
Marginal Costing and Dynamic Pricing on smart Power grids, grid topology control for
congestion mitigation, and the extension of power markets to include distribution
connected loads, generation, and resources.  He is coauthor of ``Spot Pricing of
Electricity,'' Kluwer, 1987, and has more than 100 refereed publications. His
disciplinary background is in Mathematical Economics, Optimization, and Stochastic
Dynamic Decision Making.
\end{IEEEbiography}

\end{document}